\documentclass[a4paper,twoside,11pt]{article}
\usepackage{amssymb}
\usepackage{amsmath}
\usepackage{amsthm}
\usepackage{latexsym}
\usepackage{amsfonts}
\usepackage{graphicx}
\usepackage{graphics}

\newcommand{\dis}{\displaystyle}
\theoremstyle{plain}
\newtheorem{thm}{Theorem}[section]   
\newtheorem{prop}[thm]{Proposition}

\theoremstyle{definition}
\newtheorem{rem}[thm]{Remark}

\newtheorem*{Proof}{Proof}

\newcommand{\ra}{\;\rightarrow\;}

\newcommand{\bi}{\beta}
\newcommand{\ga}{\gamma }

\newcommand{\de}{\delta }

\newcommand{\f}{\varphi}

\newcommand{\vPsi}{\varPsi}
\newcommand{\zi}{\zeta }
\newcommand{\thi}{\theta }

\newcommand{\la}{\lambda }

\newcommand{\oo}{\omega}
\newcommand{\C}{\mathbb{C}}

\newcommand{\R}{\mathbb{R}}

\newcommand{\ld}{\ldots}

\newcommand{\sm}{\smallsetminus}

\newcommand{\qb}{$\quad\blacksquare$}

\begin{document}
\pagestyle{myheadings}
\markboth{Generic versions of a result in the theory of Hardy spaces}{V. Nestoridis, E. Thirios}
\title{\bf Generic versions of a result in the theory\\ of Hardy spaces}
%
%
\author{V. Nestoridis, E. Thirios}
\date{}
\maketitle
\begin{abstract}
We show generic existence of functions $f$ in the Handy space $H^p(0<p<1)$ on the open unit disc whose primitive $F(f)$ satisfies the following.
\begin{enumerate}
\item[i)] $F(f)\in H^q$, where $q=\dfrac{p}{1-p}$.
\item[(ii)] For every $a>q$ and every $A,B\in \R$, $A<B$ it holds
\[
\sup_{0<r<1}\int^B_A|F(f)(re^{i\thi})|^ad\thi=+\infty.
\]
\item[iii)] The functions $f$ and $F(f)$ are totally unbounded and hence non-extendable.
    \end{enumerate}
Results of similar nature are valid when the space $H^p$ is replaced by localized versions of it, $0<p<1$, or intersections of such spaces.
\end{abstract}
{\bf AMS}: Classification number: 30H10 
\section{Introduction}\label{sec1}
\noindent

If there exists an object with a ``bad'' property, then a general principle is that there are many such objects and their set is big in various senses. For instance the set of irrational real numbers is uncountable while the set of rationals $Q$ is countable. Also $R\sm Q$ is bigger than $Q$ with respect to measure or category.

A classical result in the theory of Hardy spaces in the unit disc states that if $0<p<1$ and $f\in H^p$ then its primitive $F(f)$ belongs to $H^q$ where $q=\dfrac{p}{1-p}$. Furthermore, if $a>q$, then there exists $g\in H^p$ such that $F(g)\notin H^a$. According to the previously mentioned general principle we show that there are a lot of functions $g\in H^p$ with $F(g)\notin H^a$ and that their set is a $G_\de$ and dense subset of $H^p$. Furthermore we show that the function $g$ may be chosen independently of a $(a>q)$ and that generically their primitives do not belong, even not locally, to Hardy spaces of higher order $a>q$. We also obtain similar results when $H^p$, $0<p<1$, is replaced by localized versions of it, or intersections of such spaces.

In order to prove the above results we use a result of M. Siskaki which roughly speaking states that under some assumptions, if an unbounded function exists, then their set is $G_\de$ and dense in a topological vector space, not necessarily complete. Thus, our first result does not use the completeness of the Hardy spaces. However, more advanced results make use of Baire's theorem and we need the completeness of the space. Those results are results of topological genericity. We also prove a result of algebraic genericity; that is we show that the set of $g\in H^p$ $(0<p<1)$ such that for a precise $a>q=\dfrac{p}{1-p}$ the primitive $F(g)$ does not belong to $H^a$, contains a vector subspace except 0 dense in $H^p$. We do not have obtained any result of speceability; that is, that the above set contains an infinite dimensional closed vector subspace of $H^p$ except 0.

The organization of the paper is as follows. Section \ref{sec2} contains all preliminary results needed for the proofs of our results, which are contained in Section \ref{sec3}.

For the use of Baire's theorem in mathematical analysis we refer to \cite{6}, \cite{4}. For algebraic genericity and spaceability we refer to \cite{1}, \cite{2}.
\section{Preliminaries}\label{sec2}
\noindent

Let $D=\{z\in\C:|z|<1\}$ be the open unit disc. A holomorphic function $f:D\!\ra\!\C$ belongs to the Hardy space $H^p$ $(0<p<+\infty)$ if $\dis\sup_{0<r<1}\dfrac{1}{2\pi}\int^{2\pi}_0|f(re^{i\thi})|^pd\thi<+\infty$. It belongs to the Hardy space $H^\infty$ if $\dis\sup_{|z|<1|}|f(z)|<+\infty$. The space $H^\infty$ endowed with supremum norm on $D$ is a Banach space, but polynomials are not dense in this space. For $1\le p<+\infty$ the space $H^p$ endowed with the norm
\[
\|f\|_p=\sup_{0<r<1}\bigg\{\frac{1}{2\pi}\int^{2\pi}_0|f(re^{i\thi})|^pd\thi\bigg\}^{1/p}
\]
is also a Banach space. Then for $f,g\in H^p$ their distance is
\[
d_p(f,g)=\sup_{0<r<1}\bigg\{\frac{1}{2\pi}\int^{2\pi}_0|f(re^{i\thi})-g(re^{i\thi})|^p
d\thi\bigg\}^{1/p}, \ \ 1\le p<+\infty.
\]
For $0<p<1$ we endow $H^p$ with the metric
\[
d_p(f,g)=\sup_{0<r<1}\frac{1}{2\pi}\int^{2\pi}_0|f(re^{i\thi})-g(re^{i\thi})|^pd\thi, \ \ f,g\in H^p, \ \ 0<p<1
\]
and then $H^p$ becomes a topological vector space endowed with a metric invariant by translation which is complete ($F$-space). For $0<p<+\infty$ polynomials are dense in $H^p$. Also convergence in $H^p$, $0<p\le+\infty$ implies uniform convergence on each compact subset of $D$. All above are well known; see for example \cite{3}, \cite{8}.

For $a,b\in(0,+\infty]$, $a<b$ we have $H^b\subset H^a$ and the injection map is continuous. Jensen's inequality implies tha the map $a\ra\dis\sup_{0<r<1}\Big\{\dfrac{1}{2\pi}\int^{2\pi}_0|f(re^{i\thi})|^ad\thi\Big\}^{1/a}$ is increasing. Obviously we also have $\dis\sup_{0<r<1}\Big\{\dfrac{1}{2\pi}\int^{2\pi}_0|f(re^{i\thi})|^pd\thi\Big\}^{1/p}\le
\dis\sup_{|z|<1}|f(z)|$.\

Next for $0<a\le+\infty$ we consider the intersection $\bigcap\limits_{p<a}H^p$.

Convergence in this space is equivalent with convergence in all spaces $H^p$, $p<a$. Equivalently, we consider a strictly increasing sequence $p_n$ converging to $a$ and the metric in $\bigcap\limits_{p<a}H^p$ is defined by
\[
d(f,g)=\sum^\infty_{n=1}\frac{1}{2^n}\frac{d_{p_n}(f,g)}{1+d_{p_n}(f,g)}, \ \ f,g\in\bigcap_{p<a}H^p.
\]
This space is also complete, in fact an $F$-space. Obviously convergence in $\bigcap\limits_{p<a}H^p$ implies uniform convergence on each compact subset of $D$.
\begin{prop}\label{prop2.1}
Polynomials are dense in $\bigcap\limits_{p<a}H^p$, for every $a$, $0<a\le+\infty$.
\end{prop}

Fore the proof it suffices for $f\in\bigcap\limits_{p<a}H^p$ to control $d_{p_n}(f,P)$, for $n=1,\ld,N$ for any finite $N$. Because of the monotonicity of the map $p\ra\dis\sup_{0<r<1}\Big\{\dfrac{1}{2\pi}\int^{2\pi}_0|f(re^{i\thi})-P(re^{i\thi})|^pd\thi\Big\}^{1/p}$ is suffices to control $d_{p_n}(f,P)$. But this is possible, because polynomials are dense in $H^{p_N}$ since $p_N<+\infty$. Thus Proposition \ref{prop2.1} holds.

Next we present localized versions of the previous spaces (\cite{9}, \cite{10}).

Let $0<p<+\infty$ and $A,B\in\R$, $A<B$. Then a holomorphic function\linebreak $f:D\ra\C$ belongs to $H^p_{[A,B]}$ if $\dis\sup_{0<r<1}\int^B_A|f(re^{i\thi})|^p\dfrac{d\thi}{B-A}<+\infty$ and to $H^\infty_{[A,B]}$ if $\dis\sup_{0<r<1}\dis\sup_{A\le\thi\le B}|f(re^{i\thi})|<+\infty$. Because of the monotonicity of the function\linebreak $a\ra\dis\sup_{0<r<1}\Big\{\int^B_A|f(re^{i\thi})|^a\dfrac{d\thi}{B-A}\Big\}^{1/a}$ it follows $H^b_{[A,B]}\subset H^a_{[A,B]}$ for $0<a<b\le+\infty$.

Convergence in $H^p_{[A,B]}$ of a sequence $f_n$ towards $f$, $f_n,f\in H^p_{[A,B]}$ is equivalent to uniform convergence on all compact subsets of $D$ and
\[
\sup_{0<r<1}\bigg\{\int^B_A|f_n(re^{i\thi})-f(re^{i\thi})|^p\frac{d\thi}{B-A}\bigg\}^{1/p}\xrightarrow[n\ra +\infty]{}0 \ \ \text{for} \ \ 0<p<+\infty \ \ \text{and}
\]
\[
\sup_{0<r<1}\sup_{A\le\thi\le B}|f_n(re^{i\thi})-f(re^{i\thi})|\xrightarrow[n\ra +\infty]{}0 \ \ \text{for} \ \ p=+\infty.
\]
The metric giving this topology in $H^p_{[A,B]}$ is defined by
\begin{align*}
d_{p,[A,B]}(f,g)=&\sup_{0<r<1}\bigg\{\int^B_A|f_n(re^{i\thi})-g(re^{i\thi})|^p\cdot
\frac{d\thi}{B-A}\bigg\}^{1/p}\\
&+\sum^\infty_{n=2}\frac{1}{2^n}\frac{\dis\sup_{|z|\le1-\frac{1}{n}}|f(z)-g(z)|}
{1+\dis\sup_{|z|\le1-\frac{1}{n}}|f(z)-g(z)|} \ \ \text{for} \ \ 1\le p<+\infty
\end{align*}
\begin{align*}
d_{p,[A,B]}(f,g)=&\sup_{0<r<1}\int^B_A|f(re^{i\thi})-g(re^{i\thi})|^p\frac{d\thi}{B-A} \\
&+\sum^\infty_{n=2}\frac{1}{2^n}\frac{\dis\sup_{|z|\le1-\frac{1}{n}}|f(z)-g(z)|}
{1+\dis\sup_{|z|\le1-\frac{1}{n}}|f(z)-g(z)|} \ \ \text{for} \ \ 0<p<1 \ \ \text{and}
\end{align*}
\begin{align*}
d_{\infty,[A,B]}=&\sup_{0<r<1}\sup_{A\le\thi\le B}|f(z)-g(z)| \\
&+\sum^\infty_{n=2}\frac{1}{2^n}\frac{\dis\sup_{|z|\le1-\frac{1}{n}}|f(z)-g(z)|}
{1+\dis\sup_{|z|\le1-\frac{1}{n}}|f(z)-g(z)|} \ \ \text{for} \ \ p=+\infty.
\end{align*}
Obviously, convergence in $H^p_{[A,B]}$ implies uniform convergence on all compact subsets of $D$. Also for $0<a<b\le+\infty$ the injection map $H^b_{[AS,B]}\subset H^a_{[A,B]}$ is continuous. Finally these spaces are complete, in fact $F$-spaces and $H^p_{[A,B]}=H^p$ when $B-A\ge 2\pi$ and in general $H^p\subset H^p_{[A,B]}$, provided $A<B$.

Let $0<a\le+\infty$. Convergence in the space $\bigcap\limits_{p<a}H^p_{[A,B]}$ is equivalent to convergence in all $H^p_{[A,B]}$ for $p<a$. A metric in $\bigcap\limits_{p<a}H^p_{[A,B]}$ compatible with this topology is given by
\[
d(f,g)=\sum^\infty_{n=1}\frac{1}{2^n}\frac{d_{p_n,[A,B]}(f,g)}{1+d_{p_n,[A,B]}(f,g)}
\]
where $p_n$ is any strictly increasing sequence converging to $a$. This space is complete, in fact an $F$-space. Obviously convergence in $\bigcap\limits_{p<a}H^p_{[A,B]}$ implies uniform convergence on all compact subsets of $D$.

Consider the function $\dfrac{1}{(1-z)}$. It is well known that for $p>0$, this function belongs to $H^p$ if and only if $p<1$; see for example \cite{3}. The same holds for $\dfrac{1}{e^{i\oo}-z}$  for any $\oo\in\R$. If $f$ is holomorphic on the open unit disc we denote by $F(f)$ its primitive vanishing at 0; that is $F(f)(z)=\int^z_0f(\zi)d\zi$. The primitive of the function $\dfrac{1}{e^{i\oo}-z}$ is $-\log\dfrac{1}{e^{i\oo}-z}$ which belongs to $H^p$ for all $p\in(0,+\infty)$, but not to $H^\infty$. Consider the function $\dfrac{1}{(e^{i\oo}-z)^\ga}$, $\ga>0$. If $p>0$, this function belongs to $H^p$, if and only if $p<\dfrac{1}{\ga}$. For $\ga=1$ we mentioned that the primitive belongs to all $H^p$, $0<p<+\infty$ but not to $H^\infty$. For $\ga<1$ the function  $\dfrac{1}{(e^{i\oo}-z)^\ga}$ belongs to $H^1$; hence, according to Hardy's inequality \cite{3}, its primitive has an absolutely convergent Taylor series on the closed unit disc and therefore, it is in $H^\infty$. For $1<\ga$ the function $\dfrac{1}{(e^{i\oo}-z)^\ga}$, belongs to $H^\de(0<\de)$ if and only if $0<\de<\dfrac{1}{\ga}$ and the essential part of its primitive is $\dfrac{1}{(e^{i\oo}-z)^{\ga-1}}$. This primitive belongs to $H^\bi$ if and only if $\bi(\ga-1)<1$, which is equivalent to $\bi<\dfrac{1}{\ga-1}=\dfrac{\frac{1}{\ga}}{1-\frac{1}{\ga}}$.

Thus, if the function $\dfrac{1}{(e^{i\oo}-z)^\ga}$, $(1<\ga)$ belongs to an $H^\de$ (which implies $0<\de<\dfrac{1}{\ga}$), its primitive belongs to $H^{\frac{\de}{1-\de}}$ because $\dfrac{\de}{1-\de}<\dfrac{\frac{1}{\ga}}{1-\frac{1}{\ga}}$. This is a particular case of the following more general result.
\begin{thm}\label{thm2.2}
{\em(\cite{3} Theorem 5.12)}. Let $0<p<1$ and $f\in H^p$. Then its primitive $F(f)$ belongs to $H^q$, where $q=\dfrac{p}{1-p}$. Furthermore, for every $a>q$ there is a function $g$ in $H^p$ so that $F(g)\notin H^a$. Moreover, for every $W,Y,A,B\in\R$ with $W<Y$ and $A<B$ we can choose $g\in H^p\subset H^p_{[W,Y]}$ and $F(g)\notin H^a_{[A,B]}$. According to the previous discussion, it suffices to set $g(z)=\dfrac{1}{(e^{i\oo}-z)^\ga}$ with $A<\oo<B$ and $\ga\in\Big[1+\dfrac{1}{a},\dfrac{1}{p}\Big)$; the last interval is non-void because $a>q=\dfrac{p}{1-p}$ and $0<p<1$.
\end{thm}

In Section \ref{sec3} we will show that the function $g$ can be chosen independently of a $(a>q)$ and $A$ and $B$ $(A<B)$ and $W$ and $Y$ $(W<Y)$. Certainly it will not be of the form  $\dfrac{1}{(e^{i\oo}-z)^\ga}$ and it will depend on $p$, $p\in (0,1)$. Theorem \ref{thm2.2} implies easily the following.
\begin{prop}\label{prop2.3} a) Let $0<p<1$ and $f\in\bigcap\limits_{\bi<p}H^\bi$. Then $F(f)\in\bigcap\limits_{\ga<q}H^\ga$, where $q=\dfrac{p}{1-p}$. Furthermore, if $W<Y$, $A<B$, and $a>q$ there is a function $g\in\bigcap\limits_{\bi<p}H^\bi_{[W,Y]}$ such that $F(g)\notin H^q_{[A,B]}$.

b) If $f\in\bigcap\limits_{\bi<1}H^\bi$, then $F(f)\in\bigcap\limits_{\ga<+\infty}H^\ga$.

For part a) it suffices to set $g(z)=\dfrac{1}{(e^{i\oo}-z)^{1/p}}$ with $A<\oo<B$.
\end{prop}

In Section \ref{sec3} we will show that the function $g$ in Proposition \ref{prop2.3} a) can be chosen independently of $W,Y,A,B$ and $a>q$ but it will depend on $p$, $o<p<1$. In order to show that the set of these functions $g$ is big will use a version of a result of M. Siskaki (\cite{11}).
\begin{prop}\label{prop2.4}
Let $V$ be a topological vector space over the field $\R$ or $\C$. Let $X$ be a non void set and $\C^X$ the set of complex function defined on $X$.

Let $T:V\ra\C^X$ be a linear map such that for every $x\in X$ the function $V\ni f\ra T(f)(x)\in\C$ is continuous. Let $S=\{f\in V:T(f)$ is unbounded on $X\}$. Then, either $S=\emptyset$ or $S$ is a $G_\de$ and dense subset of $V$.
\end{prop}

In \cite{5} it was noticed that the previous result holds even if $T$ is not linear but it satisfies
\[
|T(f+g)(x)|\le|T(f)(x)|+|T(g)(x)| \ \ \text{andf} \ \ |T(\la f)(x)|=|\la|\,|T(f)(x)|
\]
for every scalar $\la$, every $f,g\in V$ and every $x\in X$.

The version that we will use is the following.
\begin{prop}\label{prop2.5}
Let $V$ be a topological vector space over the field $\R$ or $\C$. Let $X$ be a non empty set and $\C^X$ the set of complex functions defined on $X$. Let $T:V\ra\C^X$ be such that

1) For every $x\in X$ the function $V\ni f\ra T(f)(x)\in\C$ is continuous.

2) $|T(f+g)(x)|\le|T(f)(x)|+|T(g)(x)|$ for all $f,g\in V$ and $x\in X$.

3) For every $f\in V$ there is a function $\vPsi:\C\ra[0,+\infty)$ and a sequence $\la_n\in\C-\{0\}$, $n=1,2,\ld$ converging to 0 with $\vPsi(\la_n)\neq0$ for all $n=1,2,\ld$ and such that $|T(\la f)(x)|\ge\vPsi(\la)|T(f)(x)|$ for every $x\in X$.

We set $S=\{f\in V:T(f)$ is unbounded on $X\}$. Then, either $S=\emptyset$ or $S$ is a $G_\de$ and dense subset of $V$. In particular, the above holds if $\vPsi(\la)\equiv|\la|^c$ for some constant $c\in(0,+\infty)$ independent of $\la$.
\end{prop}
\begin{Proof}
The proof that $S$ is a $G_\de$ is omitted, because it follows simply from (1) and is similar to the proof in \cite{11}.

Suppose $S\neq\emptyset$ and let $f\in S$. Thus, $T(f)$ is unbounded on $X$. If $S$ is not dense, then there exist $g\in U$, $U$ open in $V$ with $U\cap S=\emptyset$ and $T(g)$ bounded on $X$. Therefore, there exists $M<+\infty$ such that $|T(g)(x)|<M$ for every $x\in X$.

Since $V$ is a topological vector space it holds $g+\la_nf\ra g$ as $n\ra+\infty$.

Since $U$ is open, there exists $n_0$ such that for every $n\ge n_0$ it holds $g+\la_nf\in U$. Therefore, the function $T(g+\la_{n_0}f)$ is bounded on $X$.

According to our assumptions we have
\begin{align*}
\vPsi(\la_{n_0})|T(f)(x)|&\le|T(\la_{n_0}f)(x)|=|T(g+\la_{n_0}f-g)(x)|\\
&\le|T(g+\la_{n_0}f)(x)|+|T(g)(x)|<|T(g+\la_{n_0}f)(x)|+M
\end{align*}
for all $x\in X$.

Since $\vPsi(\la_{n_0})\neq0$ it follows that $T(f)$ is bounded on $X$, which contradicts the fact that $f\in S$. Thus, $S$ is dense provided $S\neq\emptyset$. This completes the proof. \qb
\end{Proof}

We close this section mentioning that a complex function defined on the open unit disc $D$ of $\C$ is called totally unbounded (\cite{9}, \cite{10}, \cite{5}, \cite{7}) if it is unbounded on $D\cap D(\zi_0,r)$ for every $r>0$ and $\zi_0\in\partial D$ where $D(\zi_0,r)=\{z\in\C:|z-\zi_0|<r\}$. If $g$ is continuous on $D$, this is equivalent to $\dis\sup_{0<r<1}\sup_{\thi\in[A,B]}|g(re^{i\thi})|=+\infty$ for every $A,B\in\R$, $A<B$. If $g$ is holomorphic on $D$ and totally unbounded, then $g$ is non-extendable (\cite{9}, \cite{10}).
\section{The results}\label{sec3}
\noindent

It is a general principle that if there exists an object with a ``bad'' property, then most of the objects have this property. In Theorem \ref{thm2.2} we saw that for $0<p<1$ and $a>q=\dfrac{p}{1-p}$ there exists a function $g\in H^p$ such that its primitive $F(f)$ does not belong to $H^a$. We will show that the set of all such functions $g$ in $H^p$ is a $G_\de$ and dense subset of $H^p$; that is, we have topological genericity.
\begin{prop}\label{prop3.1}
Let $0<p<1$, $a>q=\dfrac{p}{1-p}$, $A<B$ and $W<Y$. Then the set $\{f\in H^p_{[W,Y]}:F(f)\notin H^a_{[A,B]}\}$ is a $G_\de$ and dense subject of $H^p_{[W,Y]}$.
\end{prop}
\noindent
{\bf Proof of Prop. \ref{prop3.1}}

Let $V=H^p_{[W,Y]}$, $X=(0,1)$ and $T:V\times X\ra\C$ be given by
\[
T(f,r)=\int^B_A|F(f)(re^{i\thi})\bigg|^a\frac{d\thi}{B-A}, \ \ f\in V, \ \ r\in X \ \ \text{if} \ \ 0<a<1
\]
and by
\[
T(f,r)=\bigg\{\int^B_A|F(f)(re^{i\thi})\bigg|^a\frac{d\thi}{B-A}\bigg\}^{1/a} \ \ f\in V, \ \ r\in X, \ \ \text{if} \ \ a\ge1.
\]
Let $0<r<1$ and $z\in\C:|z|=r$. Then $F(f)(z)=\int^z_0f(\zi)d\zi$, where the integral is over the segment $[0,z]$ which is subset of the compact disc $\{z\in\C:|z|\le r\}\subset D$. Because convergence in $H^p_{[W,Y]}$ implies uniform convergence of each compact subset of $D$, one easily can check that the assumptions of Proposition \ref{prop2.5} are verified. Furthermore, according to Theorem \ref{thm2.2} there exists $g\in V$, such that, $\dis\sup_{r\in X}|T(g,r)|=+\infty$. It follows that the set $\{f\in V:\dis\sup_{r\in X}|T(f,r)|=+\infty\}$ is a $G_\de$ and dense subset of $V$. The proof is completed.  \qb
\begin{prop}\label{prop3.2}
Let $0<p<1$, $a>q=\dfrac{p}{1-p}$ and $W<Y$. Then the set\linebreak $\{f\in H^p_{[W,Y]}$: for every $A<B$ it holds $F(f)\notin H^a_{[A,B]}\}$ is $G_\de$ and dense in $H^p_{[W,Y]}$.
\end{prop}
\begin{Proof}
The assertion that for every $A<B$, $A,B\in\R$ it holds $F(f)\notin H^a_{[A,B]}$ is equivalent that for every $A<B$, $A,B$ rational numbers it holds $F(f)\notin H^a_{[A,B]}$. According to Proposition \ref{prop3.1} for every rational numbers $A,B$, $A<B$ the set $\{f\in H^p_{[W.Y]}:F(f)\notin H^a_{[A,B]}\}$ is a $G_\de$ dense subset of the complete metric space $H^p_{[W,Y]}$. Baire's theorem implies that the denumerable intersection of these sets for all rational numbers $A,B$, $A<B$ is also a $G_\de$ and dense subset of $H^p_{[W,Y]}$. This completes the proof. \qb
\end{Proof}
\begin{prop}\label{prop3.3}
Let $0<p<1$, $q=\dfrac{p}{1-p}$, $W<Y$ and $A<B$. Then the set $\{f\in H^p_{[W,Y]}$: for every $a>q$ it holds $F(f)\notin H^a_{[A,B]}\}$ is a $G_\de$ and dense subject of $H^p_{[W,Y]}$.
\end{prop}
\begin{Proof}
Because of the monotonicity of the function
\[
a\ra\dis\sup_{0<r<1}\bigg\{\int^B_A\bigg|F(f)(re^{i\thi})\bigg|^a\frac{d\thi}{B-A}\bigg\}^{1/a} \]
the assertion that for every $a>q$ it holds $F(f)\notin H^a_{[A,B]}$ is equivalent to the assertion that for every rational number $a$, $a>q$ it holds $F(f)\notin H^a_{[A,B]}$.

Thus, the set of Proposition \ref{prop3.3} is equal to the denumerable interrection for all rational numbers $a$, $a>q$ of the sets $\{f\in H^p_{[W,Y]}:F(f)\notin H^a_{[A,B]}\}$. According to Proposition \ref{prop3.1} these sets are $G_\de$ and dense in $H^p_{[W,Y]}$. Baire's theorem implies the result. \qb
\end{Proof}
\begin{thm}\label{thm3.4}
Let $0<p<1$ and $q=\dfrac{p}{1-p}$ and $W<Y$. Then the set $\{f\in H^p_{[W,Y]}$: for every $a>q$ and every $A,B\in\R$, $A<B$ it holds $F(f)\notin H^a_{[A,B]}\}$ is a $G_\de$ and dense subset of $H^p_{[W,Y]}$.
\end{thm}
\begin{Proof}
The set of Theorem \ref{thm3.4} is equal to the intersection of the sets of Proposition \ref{prop3.3} for all rational numbers $A,B$, $A>B$. Since these sets are $G_\de$ and dense in the complete space $H^p_{[W,Y]}$ and the intersection is denumerable, Baire's theorem yields the result. \qb
\end{Proof}

If $W-Y\ge2\pi$, then $H^p_{[W,Y]}=H^p$. Combining in this case Theorem \ref{thm3.4} with Theorem \ref{thm2.2} we obtain the following.
\begin{thm}\label{thm3.5}
Let $0<p<1$ and $q=\dfrac{p}{1-p}$. Then there is a function $f$ in $H^p$ such that its primitive $F(f)$ satisfies the following.

i) $\dis\sup_{0<r<1}\int^{2\pi}_0|F(f)(re^{i\thi})|^qd\thi<+\infty$.

ii) For every $a>q$ and $A,B\in\R$, $A<B$ it holds $\dis\sup_{0<r<1}\int^B_A|F(f)(re^{i\thi})|^ad\thi=+\infty$.\\
Furthermore, the set of such functions $f$ is a $G_\de$ and dense subset of $H^p$.
\end{thm}
\begin{rem}\label{rem3.6}
In the previous results the space $H^a_{[A,B]}$ $a>q$ may be replaced by the space $H^\infty_{[A,B]}$ and similar results hold. For the proof of the analogue of Proposition \ref{prop3.1} when $a=+\infty$ it suffices to consider $V=H^p_{[W,Y]}$, $0<p<1$, $X=\{re^{i\thi}:0<r<1,A\le\thi\le B\}$ and $T(f,re^{i\thi})=F(f)(re^{i\thi})$, $f\in V$, $re^{i\thi}\in X$ and use the fact that there exists a function $g\in H^p\subset H^p_{[W,Y]}$ such that $\dis\sup_{z\in X}|F(f)(z)|=+\infty$. It suffices to consider $g(z)=\dfrac{1}{e^{i\oo}-z}$ with $A<\oo<B$. Thus, in Theorem \ref{thm3.4} when we say ``for every $a>q$'' it is the same to say ``for every $a$ in $(q,+\infty]$'' and Theorem \ref{thm3.5} has the following equivalent formulation.
\end{rem}

``Let $0<p<1$ and $q=\dfrac{p}{1-p}$. Then there is a function $f$ in $H^p$, such that its primitive $F(f)$ satisfies the following
\begin{enumerate}
\item[i)] $\dis\sup_{0<r<1}\int^{2\pi}_0|F(f)(re^{i\thi})|^qd\thi<+\infty$.
\item[ii)] For every $a\in(q,+\infty)$ and $A,B\in\R$, $A<B$ it holds $\dis\sup_{0<r<1}\int^B_A|F(f)(re^{i\thi})|^ad\thi=+\infty$ and
\item[iii)] The primitive $F(f)$ is totally unbounded on the open unit disc.
\end{enumerate}

Furthermore, the set of such functions $f$ is a $G_\de$ and dense subset of $H^p$''.

We note that for such functions $f$ and $F(f)$ are non extendable on the open unit disc.

Next we replace the space $H^p_{[W,Y]}$ by the space $\bigcap\limits_{\bi<p}H^\bi_{[W,Y]}$, $0<p\le1$. Then, Theorem \ref{thm2.2} implies that if $f\in\bigcap\limits_{\bi<p}H^\bi$ then $F(f)\in\bigcap\limits_{\ga<q}H^\ga$ with $q=\dfrac{p}{1-p}$. Let $A,B,A<B$ be given and consider the function $g=\dfrac{1}{(e^{i\oo}-z)^{1/p}}$ where $A<\oo<B$. Then $g\in\bigcap\limits_{\bi<p}H^\bi\subset\bigcap\limits_{\bi<p}H^\bi_{[W,Y}]$ and $F(g)\notin H^q_{[A,B]}$. Thus, all the previous results, extend in this case. The analogue of Theorem \ref{thm3.4} is the following.
\begin{thm}\label{thm3.7}
Let $0<p\le1$, $q=\dfrac{p}{1-p}$ and $W,Y\in\R$ with $W<Y$. Then the set $\Big\{f\in\bigcap\limits_{\bi<p}H^\bi_{[W,Y]}$: for every $A<B$ and for every $a\ge q$ it holds $F(f)\notin H^a_{[A,B]}\Big\}$ is a $G_\de$ and dense subset of the space $\bigcap\limits_{\bi<p}HYp_{[W,Y]}$.
\end{thm}
\begin{rem}\label{rem3.8}
If $W-Y\ge2\pi$, then automatically $F(f)\in\bigcap\limits_{\ga<q}H^\ga$ according to Proposition \ref{prop2.3}.
\end{rem}
\begin{rem}\label{rem3.9}
If $p=1$ and $W-Y\ge2\pi$, then Theorem \ref{thm3.6} takes the form ``The set $\Big\{f\in\bigcap\limits_{\bi<1}H^\bi:F(f)$ is totally unbounded on the open unit disc$\Big\}$ is a $G_\de$ and dense subset of the space $\bigcap\limits_{\bi<1}H^\bi$''. This is in contrast with the situation in $H^1$, where for every $f\in H^1$ the primitive $F(f)$ is bounded on the open unit disc, as it follows from Hardy's inequality (\cite{3},\cite{8}).

We close with an algebraic genericity result.
\end{rem}
\begin{thm}\label{thm3.10}
Let $0<p<1$, $a>q=\dfrac{p}{1-p}$ and $A<B$. Then the set $\{f\in H^p:F(f)\notin H^a\}$ contains a dense vector subspace of $H^p$ except zero.
\end{thm}
\begin{Proof}
Let $f_j$, $j=1,2,\ld$ be an enumeration of all polynomials with coefficients with rational real and imaginary parts. The sequence $f_j$, $j=1,2,\ld$ is dense in $H^p$ \cite3. Let $A<\oo_j<\oo_{j+1}<\dfrac{A+B}{2}$, be a sequence converging to $\dfrac{A+B}{2}$. We also assume that $\dfrac{A+B}{2}-2\pi<\oo_1$. Let $\f_j$, $j=1,2,\ld$ be a function the form $\f_j=\dfrac{1}{(e^{i\oo_j}-z)^\ga}$, where the exponent $\ga>0$ is such that $\f_j\in H^p$ and $F(\f_j)\notin H^a$.

Let $c_j\neq0$ be close enough to 0, so that $d(c_j\f_j,0)<\dfrac{1}{j}$ where $d$ is the metric in $H^p$. Then the sequence $c_j\f_j+f_j$, $j=1,2,\ld$ is dense in $H^p$, because $H^p$ does not contain isolated points.

Let $F$ be the linear span of the sequence $c_j\f_j+f_j$, $j=1,2,\ld$\;. Then $F$ is a dense vector subspace of $H^p$. It remains to show that for every non zero element $f$ of $F$ we have $F(f)\notin H^a$.

Let $f=\bi_1(c_1\f_1+f_1)+\cdots+\bi_m(c_m\f_m+f_m)$ with $\bi_m\neq0$. Let $I$ be a closed interval centered at $\oo_m$ and not containing any other $\oo_j$, $j\neq m$. Then for $j\neq m$ we have $\dis\sup_{0<r<1}\int_I|F(c_j\f_j+f_j)(re^{i\thi})|^ad\thi<+\infty$ while $\dis\sup_{0<r<1}\int_I|F(c_m\f_m+f_m)(re^{i\thi})|^ad\thi=+\infty$. It follows that $F(f)\notin H^a$. This completes the proof. \qb
\end{Proof}

In a similar way one can prove the following, because polynomials are dense in $\dis\bigcap\limits{\bi<p}H^\bi$ (see Prop \ref{prop2.1}).
\begin{thm}\label{thm3.11}
Let $0<p\le1$, $q=\dfrac{p}{1-p}$ and $A<B$. Then the set $\Big\{f\in\bigcap\limits_{\bi<p}H^\bi:F(f)\notin H^q\Big\}$ contains a dense vector subspace of $\bigcap\limits_{\bi<p}H^\bi$ except zero.
\end{thm}
\noindent
{\bf Acknowledgement} We would like to thank I. Deliyanni, T. Hatziafratis and A. Siskakis for helpful communications.
%

%
%
\noindent
V. Nestoridis \\
National and Kapodistrian University of Athens\\
Department of Mathematics\\
Panepistemiopolis, 157 84\\
Athens, Greece \medskip\\
vnestor@math.uoa.gr\\
E.Thirios@dei.com.gr

\end{document}